\documentclass[12pt]{article}
\usepackage{tikz}
\usetikzlibrary{calc}
\usepackage{bez123,calc,curves,ebezier,epic,eepic,graphicx,multiply,rotating}
\everymath{\displaystyle}
\usepackage{graphicx}
\usepackage{pst-all}
\usepackage{color}
\usepackage{amssymb}
\usepackage{amsmath}
\newtheorem{prethm}{{\bf Theorem}}

\newenvironment{thm}{\begin{prethm}{\hspace{-0.5
               em}{\bf .}}}{\end{prethm}}
\newtheorem{prelemma}{{\bf Lemma}}

\newenvironment{lemma}{\begin{prelemma}{\hspace{-0.5
               em}{\bf .}}}{\end{prelemma}}
\newtheorem{preex}{{\bf Example}}

\newtheorem{preprop}{{\bf Proposition}}

\newenvironment{prop}{\begin{preprop}{\hspace{-0.5em}{\bf .}}}{\end{preprop}}
\newtheorem{precor}{{\bf Corollary}}

\newenvironment{cor}{\begin{precor}{\hspace{-0.5
               em}{\bf .}}}{\end{precor}}
\newtheorem{preremark}{{\bf Remark}}

\newtheorem{preprob}{{\bf Problem}}

\newtheorem{predefin}{{\bf Definition}}

\newtheorem{preconj}{{\bf Conjecture}}

\newenvironment{conj}{\begin{preconj}{\hspace{-0.5
               em}{\bf . }}}{\end{preconj}}
\newtheorem{preprobb}{{\bf Problem}}

\newtheorem{prelem}{{\bf Theorem}}

\newenvironment{proof}{{\bf Proof.}\rm }{\hfill{$\Box$}}

\newtheorem{presolution}{{\bf Solution.}}

\def\newpic#1{}

\title{\vspace{-0.51cm}\Large\bf On Grundy and {\rm b}-chromatic number of some families of graphs: a comparative study}

\author{\large\bf Zoya Masih~~~~~~Manouchehr Zaker\footnote{mzaker@iasbs.ac.ir}
	\vspace{5mm}\\
	Department of Mathematics,\\
	Institute for Advanced Studies in Basic Sciences,\\
	Zanjan 45137-66731, Iran\\
}

\date{}

\begin{document}

\maketitle

\begin{abstract}
\noindent The Grundy and the {\rm b}-chromatic number of graphs are two important chromatic parameters. The Grundy number of a graph $G$, denoted by $\Gamma(G)$ is the worst case behavior of greedy (First-Fit) coloring procedure for $G$ and the {\rm b}-chromatic number ${\rm{b}}(G)$ is the maximum number of colors used in any color-dominating coloring of $G$. Because the nature of these colorings are different they have been studied widely but separately in the literature. This paper presents a comparative study of these coloring parameters. There exists a sequence $\{G_n\}_{n\geq 1}$ with limited {\rm b}-chromatic number but $\Gamma(G_n)\rightarrow \infty$. We obtain families of graphs $\mathcal{F}$ such that for some adequate function $f(.)$, $\Gamma(G)\leq f({\rm{b}}(G))$, for each graph $G$ from the family. This verifies a previous conjecture for these families.
\end{abstract}

\noindent {\bf Keywords}~~Graph coloring~.~First-Fit coloring~.~Grundy number~.~Color-dominating coloring~.~{\rm b}-chromatic number



\section{Introduction}

\noindent All graphs in this paper are undirected without any loops or multiple edges. We refer to \cite{BM} for notations and concepts not defined here. By a Grundy coloring of a graph $G$ we mean any partition of $V(G)$ into independent color classes $C_1, \ldots, C_k$ such that for each $i,j \in \{1, \ldots, k\}$ with $i<j$, each vertex in $C_j$ has a neighbor in $C_i$. The maximum integer $k$ such that there exists a Grundy coloring with $k$ colors, is called the Grundy number (also the First-Fit chromatic number) and denoted by $\Gamma(G)$ (also $\chi_{FF}(G)$). It can be observed that $\Gamma(G)$ is equal to the maximum number of colors used by the greedy coloring procedure in $G$ \cite{Z1}. The literature is full of papers concerning the extremal and algorithmic aspects of the Grundy number e.g. \cite{GL, HS, Z1, Z2}. The Grundy number is an $NP$-complete quantity even for very restricted families of graphs \cite{Z1}.
As proved in \cite{Z1}, for every integer $k$ there exists a unique tree $T_k$ (called tree atom of Grundy number $k$) such that $\Gamma(T_k)=k$ and $T_k$ is the smallest tree having Grundy number $k$. Also, for any tree $T$, $\Gamma(T)\geq k$ if and only if $T_k$ is isomorphic to a subtree of $T$. Tree atoms are also introduced in \cite{GL}.

\noindent By a color-dominating coloring of $G$ we mean any partition of $V(G)$ into independent subsets $C_1, \ldots, C_k$ such that for each $i$, the class $C_i$ contains a vertex say $v$ such that $v$ has a neighbor in any other class $C_j$, $j\not=i$. Also by a color dominating vertex $u$ of $G$  we mean the vertex $u$ has at least one neighbor in each color class of the color-dominating coloring of $G$. Denote by ${\rm{b}}(G)$ the maximum number of colors used in any color-dominating coloring of $G$. The {\rm b}-chromatic number has been widely studied in graph theory \cite{CLM, CVV, EGT, HK, IM, K2, KM, KZ, KTV}. For a survey on {\rm b}-chromatic number, we refer to \cite{JP}. In this paper we call a graph $G$ $b$-monotone if for each induced subgraph $H$ of $G$ we have ${\rm{b}}(H)\leq {\rm{b}}(G)$. This concept is similar to the concept of quasi-monotonous graphs, introduced in \cite{K2}, where a graph $G$ is called quasi-monotonous if for any two subgraphs $H_1$ and $H_2$ of $G$ such that $H_1$ is a subgraph of $H_2$, one has ${\rm{b}}(H_1)\leq {\rm{b}}(H_2)$.

\noindent Obviously, every quasi-monotonous graph is also $b$-monotone. There is a useful quantity denoted by $m(G)$ which is used in study of {\rm b}-chromatic number. Let the vertex degrees of $G$ be ordered in a non-increasing form $d_1 \geq d_2 \geq \ldots \geq d_n$. Define $m(G)=\max \{k: d_k\geq k-1\}$.
It is easily seen that ${\rm{b}}(G)\leq m(G)$, since otherwise if ${\rm{b}}(G)>m(G)$, there is a color $c$ such that all the vertices with color $c$ have degree $\leq m-1$. Hence all vertices of color $c$ have degree less that $\rm{b}(G)-1$ and none of them can be a color dominating vertex, a contradiction.

\noindent A natural question concerning comparison of Grundy and {\rm b}-chromatic numbers is to explore or generate families of graphs $\{G_n\}_{n\geq 1}$ and $\{H_n\}_{n\geq 1}$ such that ${\rm{b}}(G_n)-\Gamma(G_n)\rightarrow \infty$ and $\Gamma(H_n)-{\rm{b}}(H_n)\rightarrow \infty$. It was proved in \cite{Z3} that both of the above-mentioned situations happen in the universe of graphs. Given any family $\mathcal{F}$, it was also proved in \cite{Z3} that there exits a function $f(.)$ such that $\Gamma(G)\leq f({\rm{b}}(G))$, for each graph $G\in \mathcal{F}$ if and only if for any sequence $\{G_n\}_{n\geq 1}$ from $\mathcal{F}$, $\Gamma(G_n)\rightarrow \infty$ implies ${\rm{b}}(G_n)\rightarrow \infty$. The following conjecture was made in \cite{Z3}.

\begin{conj}
There exits a function $f(.)$ such that $\Gamma(G)\leq f({\rm b}(G))$ for any $b$-monotone graph $G$.\label{conj}
\end{conj}

\noindent As proved in \cite{Z3}, the conjecture is not valid for the family of non $b$-monotone graphs. In this paper we prove validity of the conjecture for some families of graphs. But in the following we first obtain a sequence of graphs with limited {\rm b}-chromatic number and unbounded Grundy number.

\noindent Let $m\geq 2$ be any fixed integer and $n$ an arbitrary positive integer. Consider first the complete $m$-partite graph $K_{n,\ldots, n}$ in which each partite set has exactly $n$ vertices. Name the partite sets by $B_1, \ldots, B_m$. Set $B_i=\{v_{i,1}, v_{i,2}, \ldots, v_{i,n}\}$. Consider $n-1$ (edge disjoint) cliques $A_1, \ldots, A_{n-1}$, where $V(A_{j}) = \{v_{1,j}, v_{2,j}, \ldots, v_{m,j}\}$. Now, remove the edges of these cliques from the graph and call the resulting graph $G_{m,n}$. In the following we prove that $G_{m,n}$ satisfies the desired properties.

\begin{prop}
For any $m\geq 2$ and $n\geq 1$, $\Gamma(G_{m,n}) \geq m+(n-1)$ and ${\rm{b}}(G_{m,n})=m$.\label{Gmn}
\end{prop}

\noindent\begin{proof}
Assign the colors $1, 2, \ldots, n-1$ to the vertices $v_{i,1}, v_{i,2}, \ldots, v_{i,n-1}$, respectively, in each partite set $B_i$, $i\in \{1, \ldots, m\}$. For each $i$, assign the color $n+i-1$ to $v_{i,n}$. The resulting coloring is a Grundy coloring using $n+m-1$ colors. It follows that $\Gamma(G_{m,n}) \geq n+m-1$.

\noindent To prove the equality concerning ${\rm{b}}$-chromatic number, first assign color $j$ to all vertices in $B_j$, $j\in \{1, 2, \ldots, m\}$. In this coloring, the vertices $v_{1,n}, v_{2,n}, \ldots, v_{m,n}$ are color-dominating vertices with $m$ different colors. It implies ${\rm{b}}(G)\geq m$. Assume on the contrary that ${\rm{b}}(G)=m+t$, for some $t\geq 1$. Let $C$ be a ${\rm b}$-coloring of the graph with $m+t$ colors. Assume that the colors in the clique $A_n$ are $i_1, i_2, \ldots, i_m$. At least one color, say $c$, $c \in \{1, \ldots, m+1 \}$, is missing in $A_n$. It implies that no color-dominating vertex of color $c$ exists, because if a vertex $v_{j,k}, j \in \{1, \ldots, m \}$ and $k \in \{1, \ldots, n-1 \}$, is color-dominating of color $c$, then this vertex has a neighbor of color $c(v_{j,n})$, by the definition of ${\rm b}$-coloring. But the only partite set containing a vertex of color $c(v_{j,n})$ is $B_j$, a contradiction.	
\end{proof}

\noindent We now show that $G_{m,n}$ is not $b$-monotone.
\begin{prop} For any $m\geq 2$ and $n\geq 1$, $G_{m,n}$ is not $b$-monotone.
\end{prop}

\noindent \begin{proof}
By Proposition \ref{Gmn}, ${\rm{b}}(G_{m,n})=m$. Then, to prove Proposition 2 we obtain a subgraph of $G_{m,n}$ with {\rm b}-chromatic number more than $m$. Remove the vertices $v_{m,3}, v_{m,4}, \ldots, v_{m,n}$ from the partite set $B_m$ and call the resulting subgraph $G'$. Note that after this removal, $G'$ has only two vertices in the $m$-th partite set. Now, in each partite set $B_i$, $1 \leq i \leq m-1$, assign color $i$ to all vertices $v_{i,3}, v_{i,4}, \ldots, v_{i,n}$ and color the vertices $v_{i,1}, v_{i,2}$, respectively with colors $m, m+1$ in each partite set $B_i$, $1 \leq i \leq m$. This is a ${\rm b}$-coloring of $G'$ with $m+1$ colors and the vertices $v_{1,n}, v_{2,n}, \ldots, v_{m-1,n}, v_{m,1}, v_{m,2}$ are color dominating vertices.
\end{proof}

\noindent Let $G$ be any graph and $C$ a Grundy coloring of $G$ using $k$ colors. We define a colored subgraph of $G$ as follows. Let $v_k$ be a vertex of color $k$ in $C$. Choose a set consisting of $k-1$ neighbors of $v_k$, say $v_1, \ldots, v_{k-1}$ with distinct colors $1, \ldots, k-1$, respectively.
Define $L_1=\{v_k\}$ and $L_2=\{v_1, \ldots, v_{k-1}\}$. Now, for each $i$ and $j$ with $1\leq i < j \leq k-1$, any vertex of color $j$ in $L_2$ needs a neighbor of color $i$ in $C$. If such a vertex say $u$ is not found in $L_2$ then put $u$ in a newly defined set $L_3$. The set $L_3$ consists only of such vertices. Repeat the above procedure for $L_3$. Any vertex $v\in L_3$ of color say $j$ needs a neighbor of color $i$ for each $i<j$. If such a neighbor $w$ is not found in $L_2\cup L_3$ then put $w$ in a newly defined set $L_4$. The set $L_4$ consists only of such vertices. We continue this procedure and obtain the other sets $L_5, \ldots$.
Let $t$ be an integer such that in $L_t$ the corresponding set $S$ is empty.
For each $i\in \{1, \ldots, t\}$, define $H_i$ as the subgraph of $G$ induced by $L_1\cup \ldots \cup L_i$. Set also $H=G[L_1\cup \ldots \cup L_t]$. We call $H$ a subgraph corresponding to the Grundy coloring $C$. For each $i\geq 1$, the presence of a vertex $v$ of color $j$ in $L_i$ means that there exists a vertex $u$ in $L_{i-1}$ such that $v$ is adjacent to $u$ and the color of $u$ is greater than $j$. In this situation $v$ is said to be a child of $u$ and we write $v\in CH(u)$. Clearly, $CH(u)\subseteq L_i$ for every $u\in L_{i-1}$.

\noindent {\bf The outline of the paper is as follows.} In this paper we prove Conjecture \ref{conj} for trees, cactus graphs, $(K_4 \setminus e, C_4)$-free {\rm b}-monotone graphs and graphs of girth at least $6$. In Sections $2$ and $3$, we obtain an upper bound in terms of {\rm b}-chromatic number for the Grundy number of trees and cacti, respectively. These bounds are almost sharp. In Section $4$, we obtain a similar result for graphs of girth five and six and for another family which is defined by some forbidden induced subgraphs.

\section{Results for trees}

\noindent It was proved in \cite{IM} that $m(T)-1 \leq {\rm{b}}(T) \leq m(T)$ for any tree $T$. A vertex $v$ of $T$ is called dense if $d_T(v)\geq m(T)-1$. A tree $T$ is said to be pivoted in \cite{IM} if $T$ has exactly $m(T)$ dense vertices, and $T$ contains a distinguished vertex $v$ such that:\\
\noindent $(i)$ $v$ is not dense.\\
\noindent $(ii)$ Each dense vertex is adjacent either to $v$ or to a dense vertex adjacent to $v$.\\
\noindent $(iii)$ Any dense vertex adjacent to $v$ and to another dense vertex has degree $m(T)-1$.

\noindent Note that a corollary of $(iii)$ is the following.

\noindent $(iii)'$ The vertex $v$ must be adjacent to at least two dense vertices.

\noindent The following result was proved in \cite{IM}.

\begin{prop}
For any non-pivoted tree $T$, ${\rm{b}}(T)=m(T)$.\label{tree-pivot}
\end{prop}

\begin{figure}
\begin{center}
\begin{tikzpicture}
\draw[black, thick] (0,0)-- (1,0);
\draw[black, thin] (3,0)-- (4,0);
\draw[black, thin] (6,0)-- (7,0);
\draw[black, thick] (9,0)-- (10,0);
\draw[black, thin] (10,0)-- (11,0);
\draw[black, thick] (11,0)-- (12,0);
\draw[black, thick] (3,0)-- (3,1);
\draw[black, thick] (4,0)-- (4,1);
\draw[black, thin] (6,0)-- (6,1);
\draw[black, thin] (7,0)-- (7,1);
\draw[black, thick] (6,0)-- (6,-1);
\draw[black, thick] (7,0)-- (7,-1);
\draw[black, thick] (6,1)-- (6,2);
\draw[black, thick] (7,1)-- (7,2);
\draw[black, thick] (9,0)-- (10,0);
\draw[black, thick] (9,1)-- (10,1);
\draw[black, thick] (9,2)-- (10,2);
\draw[black, thin] (10,1)-- (10,2);
\draw[black, thin] (10,0)-- (10,1);
\draw[black, thin] (10,-1)-- (10,0);
\draw[black, thick] (9,2)-- (10,2);
\draw[black, thick] (9,1)-- (10,1);
\draw[black, thin] (11,-1)-- (11,0);
\draw[black, thick] (9,-1)-- (10,-1);
\draw[black, thick] (11,-1)-- (12,-1);
\draw[black, thin] (11,0)-- (11,1);
\draw[black, thick] (11,1)-- (12,1);
\draw[black, thin] (11,1)-- (11,2);
\draw[black, thick] (11,2)-- (12,2);
\filldraw [black] (1,0) circle(2pt)
node [anchor=west]{};
\filldraw [gray] (0,0) circle(2pt)
node [anchor=south]{};
\filldraw [gray] (3,0) circle(2pt)
node [anchor=south]{};
\filldraw [gray] (4,0) circle(2pt)
node [anchor=south]{};
\filldraw [black] (4,1) circle(2pt)
node [anchor=south]{};
\filldraw [gray] (7,0) circle(2pt)
node [anchor=south]{};
\filldraw [gray] (7,1) circle(2pt)
node [anchor=south]{};
\filldraw [gray] (6,1) circle(2pt)
node [anchor=south]{};
\filldraw [black] (3,1) circle(2pt)
node [anchor=south]{};
\filldraw [gray] (6,0) circle(2pt)
node [anchor=south]{};
\filldraw [black] (7,2) circle(2pt)
node [anchor=south]{};
\filldraw [black] (6,2) circle(2pt)
node [anchor=south]{};
\filldraw [black] (7,-1) circle(2pt)
node [anchor=south]{};
\filldraw [black] (6,-1) circle(2pt)
node [anchor=south]{};
\filldraw [gray] (10,0) circle(2pt)
node [anchor=south]{};
\filldraw [gray] (11,0) circle(2pt)
node [anchor=south]{};
\filldraw [gray] (10,1) circle(2pt)
node [anchor=south]{};
\filldraw [gray] (11,1) circle(2pt)
node [anchor=south]{};
\filldraw [gray] (10,2) circle(2pt)
node [anchor=south]{};
\filldraw [gray] (11,2) circle(2pt)
node [anchor=south]{};
\filldraw [gray] (10,-1) circle(2pt)
node [anchor=south]{};
\filldraw [gray] (11,-1) circle(2pt)
node [anchor=south]{};
\filldraw [black] (9,-1) circle(2pt)
node [anchor=south]{};
\filldraw [black] (9,0) circle(2pt)
node [anchor=south]{};
\filldraw [black] (9,1) circle(2pt)
node [anchor=south]{};
\filldraw [black] (9,2) circle(2pt)
node [anchor=south]{};
\filldraw [black] (12,-1) circle(2pt)
node [anchor=south]{};
\filldraw [black] (12,0) circle(2pt)
node [anchor=south]{};
\filldraw [black] (12,1) circle(2pt)
node [anchor=south]{};
\filldraw [black] (12,2) circle(2pt)
node [anchor=south]{};

\end{tikzpicture}
\caption{$T_2, T_3, T_4, T_5$, from left to right}\label{Ti}
\end{center}
\end{figure}
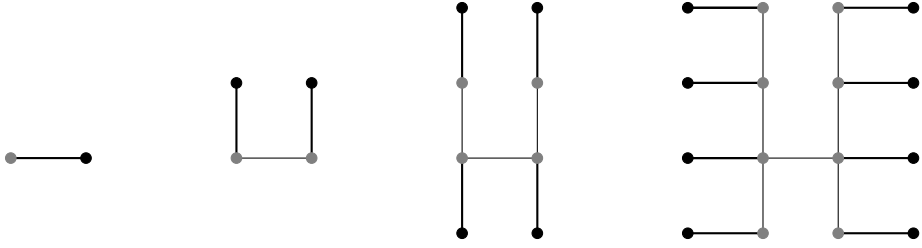

\noindent We are going to prove Proposition \ref{bound-atom} and then Corollary \ref{bound-tree} concerning all trees. We need a proposition concerning $m(T_k)$ of tree atoms $T_k$. We define the $k$-atom tree $T_k$, by induction on $k$. The $T_1$ and $T_2$ are isomorphic to the complete graphs on one and two vertices, respectively. Suppose that we have constructed $T_k$ on $n$ vertices. To construct $T_{k+1}$, attach a leaf to each vertex of $T_k$. It follows that $|V(T_k)|=2^{k-1}$. The tree $T_k$ is illustrated in Figure \ref{Ti}, for $k=2, 3, 4, 5$.
It can be shown that for each $k$, $\Gamma(T_k)=k$ and for any tree $T$, $\Gamma(T)\geq k$ if and only if $T_k$ is isomorphic to a subtree of $T$. Also note that for every positive integer $k$, there exists a unique non-negative integer $i$ such that $2^i+i\leq k \leq 2^{i+1}+i$. The reason is that the set of natural numbers $\Bbb{N}$ is partitioned into disjoint integer intervals ${\bigcup}_{i\geq 0} [2^i+i, 2^{i+1}+i]$. We need the following proposition.

\begin{prop}
Let $k$ be any positive integer and let $i$ be the unique integer such that $2^i+i \leq k \leq 2^{i+1}+i$. Then $m(T_k)=k-i$.\label{atom-m}
\end{prop}

\noindent \begin{proof}
Assume that the vertex degrees in $T_k$ are in a non-increasing form $d_1\geq d_2\geq \cdots \geq d_n$. For $T_2$ and $T_3$ the related lists are $1, 1$ and $2, 2, 1, 1$, respectively. Note that if $d_1, d_2,  \ldots, d_n$ is a degree sequence of $T_k$, then the degree sequence of $T_{k+1}$ is $d_1+1, d_2+1, \ldots, d_n+1, \underbrace{1, \ldots, 1}_n$, by its construction.

It follows that the size of the degree sequence of $T_k$ is $ 2^{k-1}$, $|T_{k}|= 2^{k-1}$,
and the degree sequence is as follows:
$$\underbrace{k-1, k-1}_2, \underbrace{k-2, k-2}_2, \underbrace{k-3,\ldots, k-3}_4, \ldots, \underbrace{2,\cdots, 2}_{2^{k-3}},\underbrace{ 1, \ldots,1}_{2^{k-2}}$$

\noindent We obtain the following result for the non-increasing degree sequence $d_1, d_2, \ldots, d_n$ of $T_k$. Let $j$ be any integer with $1< j \leq n$ such that $2^{t-1} < j \leq 2^t$, for some integer $t$. Then $d_j=k-t=k-\lceil \log{j} \rceil$ also we have  $d_1=k-1$. We show this by the induction on $k$. At first we consider $d_1$. In $T_1$, which is a single vertex, we have $d_1=0$. Suppose that in $T_i$ we have $d_1={i-1}$. As we mentioned before, in the degree sequence of $T_{i+1}$, the greatest degree is $d_{1}+1$ which means in $T_{i+1}$ we have $d_1=i$.
Consider now $d_j$ for $1< j \leq n$.
For $k=2$, $T_2$, we have $d_2=1=k-\lceil \log{2} \rceil$. In $T_i$, we have $|T_{i}|= 2^{i-1} $. Assume that the vertex degrees in $T_i$ are in a non-increasing form $d_1, d_2, \cdots, d_{2^{i-1}}$.  Suppose that for any integer $j$, with $1< j \leq 2^{i-1}$, we have $d_j=i-\lceil \log{j} \rceil$. In $T_{i+1}$ the vertex degree sequence will be as $d_{1}+1, d_{2}+1, \cdots, d_{2^{i-1}}+1, 1, \cdots, 1$. Denote by $d'_j$ the $j$-th vertex degree in the degree sequence of $T_{i+1}$. For any $j$ with $ 2\leq j \leq 2^{i-1}$ we have $d'_j=d_j +1=i-\lceil \log j \rceil+1$.
Also for  any $j$ with $ 2^{i-1}< j \leq 2^{i}$ we have $d'_j=1$ which means $d'_j=i-\lceil \log j \rceil+1$.

\noindent Write for simplicity $k-i=p$. To prove the proposition, it suffices to show $d_p\geq p-1$ and $d_{p+1} < p$. Let $k= 2^{i}+i+r$ for some $0\leq r \leq 2^{i+1}-2^{i}$, hence $p=k-i= 2^{i}+r$. By the above expression for $d_j$ we obtain $d_p=d_{2^{i}+r}=k-i-1=p-1$. From the other side,
$$d_{p+1}=d_{2^i+r+1}=k-\lceil \log{(2^i+r+1)}\rceil=(2^i+i+r)-\lceil \log{(2^i+r+1)}\rceil <2^i+r=p.$$
\end{proof}

\noindent We use Proposition \ref{atom-m} in the following result.

\begin{prop}
For any positive integer $k$, $k-\lfloor \log (k-1) \rfloor \leq {\rm{b}}(T_k)$.\label{bound-atom}
\end{prop}

\noindent \begin{proof}
\noindent It follows from Proposition \ref{atom-m} that
$k-\lfloor \log (k-1) \rfloor \leq m(T_k) \leq k- \lfloor \log k \rfloor +1$ since $\lfloor \log k \rfloor -1 \leq i \leq \lfloor \log(k-1) \rfloor$.
Note that dense vertices of $T_k$ induces a subtree of $T_k$. Hence every other vertex of $T_k$ has at most one neighbor in this subtree. Therefore the distinguished vertex $v$ does not exist by the condition $(iii)'$ and $T_k$ is not pivoted.
By Proposition \ref{tree-pivot}, ${\rm{b}}(T_k)=m(T_k)$ and the above-mentioned inequalities for $m(T_k)$ hold for ${\rm{b}}(T_k)$, i.e.
\begin{center}
$k-\lfloor \log (k-1) \rfloor \leq {\rm{b}}(T_{k}) \leq k- \lfloor \log k \rfloor +1$.
\end{center}
\end{proof}

\noindent It can be easily proved that trees are $b$-monotone. This result can also be deduced from a result of \cite{K2} that graphs of girth at least five are quasi-monotonous. We obtain the following corollary which validates Conjecture \ref{conj} for trees.

\begin{cor}
Let $T$ be a tree graph. Then $$\Gamma(T)-\lfloor \log (\Gamma(T)-1) \rfloor \leq {\rm{b}}(T).$$\label{bound-tree}
\end{cor}	

\section{Cactus graphs}

\begin{figure}
\begin{center}
\begin{tikzpicture}
\draw[black, thick] (0,0)-- (-1,0);
\draw[black, thick] (0,0)-- (1,-1);
\draw[black, thick] (0,0)-- (1,1);
\draw[black, thick] (-1,0)-- (-2,1);
\draw[black, thick] (-1,0)-- (-2,-1);
\draw[black, thick] (0,-3)-- (1,-2);
\draw[black, thick] (0,-3)-- (1,-4);
\draw[black, thick] (0,-3)-- (-1,-3);
\draw[black, thick] (-1,-3)-- (-2,-2);
\draw[black, thick] (-1,-3)-- (-2,-4);
\draw[black, thick] (2,-1.5)-- (1,-2);
\draw[black, thick] (2,-1.5)-- (1,1);
\draw[black, thick] (2,-1.5)-- (1,-1);
\draw[black, thick] (2,-1.5)-- (1,-4);
\filldraw [black] (2,-1.5) circle(2pt)
node [anchor=west]{$ v$};
\filldraw [black] (1,-2) circle(2pt)
node [anchor=south]{};
\filldraw [black] (1,-4) circle(2pt)
node [anchor=south]{};
\filldraw [black] (-2,-4) circle(2pt)
node [anchor=south]{};
\filldraw [black] (-2,-2) circle(2pt)
node [anchor=south]{};
\filldraw [black] (-1,-3) circle(2pt)
node [anchor=south]{};
\filldraw [black] (0,-3) circle(2pt)
node [anchor=south]{};
\filldraw [black] (-2,-1) circle(2pt)
node [anchor=south]{};
\filldraw [black] (-2,1) circle(2pt)
node [anchor=south]{};
\filldraw [black] (-1,0) circle(2pt)
node [anchor=south]{};
\filldraw [black] (1,-1) circle(2pt)
node [anchor=south]{};
\filldraw [black] (1,1) circle(2pt)
node [anchor=south]{};
\filldraw [black] (0,0) circle(2pt)
node [anchor=south]{};

\end{tikzpicture}
\caption{A cactus graph which is not {\rm b}-monotone.}\label{noncuctus}
\end{center}
\end{figure}
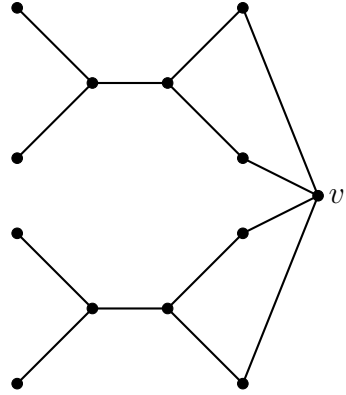

\noindent A graph $G$ is cactus if each block of $G$ is isomorphic to an edge or a cycle, where by a block we mean any maximal $2$-connected subgraph in the graph. The {\rm b}-chromatic number of cactus graphs was studied in \cite{CLM}. Cactus graphs are not necessarily $b$-monotone. Figure \ref{noncuctus} illustrates a cactus graph which is not $b$-monotone, where ${\rm b}(G\setminus v)=4$ and
${\rm b}(G)=3$.

\noindent Suppose that $G$ is any cactus graph with $\Gamma(G)=k$ and $C$ be any Grundy coloring of $G$ using $k$ colors. Let $H$ be the subgraph of $G$ corresponding to $C$, as defined previously in the introduction.

\begin{lemma}
Let $G$ be a cactus graph and $i$ with $2\leq i \leq t$ be a fixed integer and $v$ be any vertex of $L_i$. Denote by $CH(v)$ the set of children of $v$. Let $w, z\in CH(v)$. If $w$ is adjacent to some vertex in $H_{i+1}\setminus CH(v)$, then $z$ is not adjacent to any vertex of $H_{i+1}\setminus CH(v)$.\label{H-cactus}
\end{lemma}

\noindent \begin{proof}
Suppose that $x$ and $y$ are two vertices in $H_{i+1}\setminus CH(v)$ such that $w$ is adjacent to $x$ and $z$ is adjacent to $y$. There exist three paths $P$, $Q$ and $R$ from $v_k$ to $v$, $x$ and $y$, respectively. We obtain two cycles $C: v_kQxwvPv_k$ and $C': v_kRyzvPv_k$ which intersect in at least two vertices. In other possibilities too, we obtain two cycles intersecting in at least two vertices. This contradicts the fact that $G$ is a cactus graph.
\end{proof}

\noindent The next theorem validates Conjecture \ref{conj} for the cactus graphs.

\begin{thm}
Let $G$ be a $b$-monotone cactus graph. Set $\Gamma(G)=k$. Then ${\rm{b}}(G)\geq k- 2\lfloor \log k\rfloor$.\label{bound-cactus}
\end{thm}

\noindent \begin{proof}
The assertion obviously holds for $k\leq 6$. For $7\leq k \leq 9$ we can easily check that a cactus graph of Grundy number $k$ has either a triangle or an induced $P_5$ or $C_5$. Each of which has {\rm b}-chromatic number $3$. Since $G$ is $b$-monotone then the assertion holds for $k\leq 9$.

\noindent Assume hereafter that $k>9$. Let $C$ be any Grundy coloring of $G$ with $k$ colors and $H$ be the subgraph corresponding to $C$. Recall that $H$ is consisted of the levels $L_1, \ldots, L_t$, for some $t$ in which $L_1=\{v_k\}$ and $L_2=\{v_1, \ldots, v_{k-1}\}$. The color of $v_i$ is $i$ and $v_i$ is adjacent to $v_k$ for each $i\not=k$. Let $p = k- 2\lfloor \log k\rfloor$. In the following we obtain an induced subgraph of $G$ which admits a {\rm b}-coloring using $p$ colors. In fact we obtain such an induced subgraph of $G$ by recoloring some portion of $G$ so that the resulting coloring is a {\rm b}-coloring with $p$ colors.

\noindent \textbf{Recoloring Process:}

\noindent As we mentioned, we are going to perform a recoloring on a portion of $G$. The colors are coming from a set of colors, $S$. We describe the recoloring process by using cases. Each case corresponds to the recoloring process for $L_1, L_2, L_3, L_4$ and $L_i, i \geq 5$. Then we prove that for each $i$ with $p - 1 \leq i \leq k - 1$, the vertex $v_i$ has all colors $1, \ldots, p$ in its closed neighborhood $N[v_i]$ and after recoloring of $L_4$, any recolored vertex $u$ of $L_3$ with initial color $p - 1$ has all colors $1, \ldots, p$ in its closed neighborhood. Likewise, we prove that  after recoloring of $L_i$, each recolored vertex $u$ of $L_{i-1}$ with $c(u)\geq p-1$ has all colors $1, \ldots, p$ in its closed neighborhood.

\noindent Let $S=\{1, \ldots, 2p-k-1, \ldots, p-1, p\}$ be an ordered set of colors to be used in the recoloring process. At each step of the recoloring process, some colors are removed from $S$ and $S$ is updated by this removal. Note that the procedure of recoloring stops whenever $S$ is empty. Also note that the following recoloring process may cause two adjacent vertices with a same new color. This problem is resolved later by a trick.  Denote the resulting recoloring by $C'$ and the new color of any vertex $v$ is denoted by $c'(v)$.
 In $L_1$, assign the color $p$ to $v_k$ and
remove the value $p$ from $S$. In $L_2$, recolor $v_{p-1}, \ldots ,v_{k-1}$ by $2p-k-1, \ldots, p-1$, respectively.  Remove the values of the set  $\{2p-k-1, \ldots, p-1\}$ from $S$.
 The color of other vertices in $L_2$ remains unchanged. The situation is depicted in Figure \ref{fig2}.
We should mention that for all values $13<k$, we have $ 2p-k-1=k-4\lfloor \log k \rfloor-1 \geq 1$. For the case $10\leq k \leq13$ (in which $ 2p-k-1=k-4\lfloor \log k \rfloor-1< 1$), we have $S=\{1, \ldots, p-1, p\}$ and we recolor vertices $v_7, v_8, \ldots ,v_{k-1},v_k$ by $1, 2, \ldots, p-1, p$, respectively. It means that in the vertex $v_7 \in L_2$ the set $S$ becomes empty and therefore the procedure is finished in $ L_2$ .

\noindent We recolor the vertices of $L_3$. Assume that the vertices of $L_3$ are presented according to an arbitrary but fixed ordering. Let $w$ be a vertex of $L_3$ whose color is $c(w)$ in the Grundy coloring $C$. Then $w$ belongs to $CH(v_i)$, for some $v_i\in L_2$. Assume that $c(w) \leq p-2$. If $c(w)$ is the same as the new color of $v_i$, then change the color $c(w)$ to $p-1$.
Otherwise, the color of $w$ remains unchanged. Note that if $c(w)\geq p-1$ and the new color of $v_i$ equals $c(w)$ then $c(w)=p-1$.
In case that none of the above situations happen, assign the greatest color in $S$ to $w$, update $S$ by removing this color from $S$ and go to another vertex in $L_3$. In this step for each $i$ with $p-1 \leq i \leq k$, the vertex $v_i$  has all colors $1, \ldots, p$ in its closed neighborhood. This is because, for each $i$ with $p-1 \leq i \leq k-1$, the vertex $v_i$  has all colors $1, \ldots, p-2$ in its neighborhood in the primary Grundy coloring.
In the recoloring process these colors remain unchanged unless for the neighbor of $v_i$, say $u$, with color $c(u)=c'(v_i)$. Here in the recoloring process the vertex $u$ receives new color $p-1$, i.e. $c'(u)=p-1$. Also $v_k$ with new color $p$ is in the neighborhood of $v_i$.  Therefore in this step $v_i$ has all colors $1, \ldots, p$ in its closed neighborhood, $N[v_i]$.
In the case $i=k$, $v_k$ is adjacent to $v_1, \ldots, v_{p-2}$ with colors respectively $1, \ldots, p-2$ and note that these colors do not change in the recoloring process. Also $v_k$, itself, has new color $p$ and is adjacent to $v_{k-1}$ which is recolored newly by $p-1$. Therefore in this step $v_k$ has all colors $1, \ldots, p$ in its closed neighborhood, $N[v_k]$.

\noindent Now consider the vertices of $L_4$ according to a fixed ordering. Let $y$ be a vertex in $L_4$ with color $c(y)$ in the coloring $C$. Then, $y$ belongs to $CH(x)$ for some $x \in L_3$ and $x$ belongs to $CH(v_i)$ for some $v_i \in L_2$. If the color of $x$ is unchanged, do not change the color of $y$.
Now suppose that $x$ has been recolored. If $c(y) \leq p-2$ and $c'(x)\neq c(y)$ or $c'(v_i) \neq c(y)$ then, the color of $y$ remains unchanged. If $c(y)\leq p-2$ and $x$ is recolored by $c(y)$ then, we recolor $y$ by $p$.
If $c(y) \leq p-2$ and $v_i$ is recolored by $c(y)$ then, we recolor $y$ by $p-1$. If none of the above situations happens (i.e. $c(y)\geq p-1$) then, assign the greatest color in $S$ to $y$ as its new color. Update $S$ by removing the assigned color from it and go to another vertex in $L_4$.
Note that after recoloring of $L_4$, any recolored vertex $u$ of $L_3$ with $c(u)\geq p-1$ has all colors $1, \ldots, p$ in its closed neighborhood. We argue this for the case $c(u)=p-1$.
Let $u\in CH(v_i)$. The vertex $u$  has all colors $1, \ldots, p-2$ in its neighborhood in the primary Grundy coloring. In the recoloring process these colors remain unchanged unless for some neighbors of $u$, say $y$ and $z$, for which we have $c(y)=c'(u)$ and  $c(z)=c'(v_i)$. Here in the recoloring process the vertices $y$ and $z$ receive new colors $p$ and $p-1$, i.e.  $c'(y)=p$ and $c'(z)=p-1$.

\begin{figure}
\label{fig2}
\begin{center}
\begin{tikzpicture}
\draw[black, thick] (-3,0)-- (-6,-1);
\draw[black, thick] (-3,0)-- (-1,-1);
\draw[black, thick] (-3,0)-- (-4.5,-1);
\draw[black, thick] (-3,0)-- (-2,-1);

\draw[black, dashed] (-2,-1)-- (-1,-1);
\draw[black, dashed] (-6,-1)-- (-4.5,-1);
\filldraw [black] (-3,0) circle(2pt)
node [anchor=south]{$v_{k}$};

\filldraw [black] (-1,-1) circle(2pt)
node [anchor=west]{$v_{k-1}$};

\filldraw [black] (-2,-1) circle(2pt)
node [anchor=north]{$v_{p-1}$};

\filldraw [black] (-4.5,-1) circle(2pt)
node [anchor=north]{$v_{p-2}$};

\filldraw [black] (-6,-1) circle(2pt)
node [anchor=east]{$v_1$};


\draw[->,line width=2pt,black]%
(-0.5,0) -- (0.5,0) ;

\draw[black, thick] (4,0)-- (1,-1);
\draw[black, thick] (4,0)-- (5,-1);
\draw[black, thick] (4,0)-- (2.5,-1);
\draw[black, thick] (4,0)-- (6,-1);

\filldraw [black] (4,0) circle(2pt)
node [anchor=south]{$p$};
\filldraw [black] (5,-1) circle(2pt)
node [anchor=north]{${2p-k-1}$};
\filldraw [black] (1,-1) circle(2pt)
node [anchor=east]{$1$};

\filldraw [black] (2.5,-1) circle(2pt)
node [anchor=north]{$p-2$};

\filldraw [black] (6,-1) circle(2pt)
node [anchor=west]{$p-1$};

\draw[black, dashed] (5,-1)-- (6,-1);
\draw[black, dashed] (1,-1)-- (2.5,-1);

\end{tikzpicture}
\end{center}
\caption{Recoloring process in cactus graphs (related to Theorem \ref{bound-cactus})}
\end{figure}
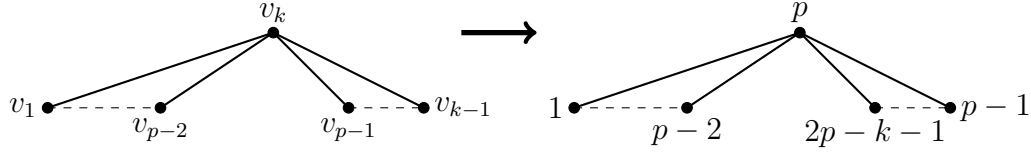
	
\noindent Now we explain a general case $L_i$. Let $y$ be any vertex in $L_i$. Then $y$ belongs to $CH(x)$ for some vertex $x \in L_{i-1}$ and $x$ belongs to $CH(w)$ for some $w \in L_{i-2}$. If $x$ has not been recolored in previous levels, then keep the color of $y$ unchanged. Suppose that the color of $x$ has been changed in the previous stage. If $c(y) \leq p-2$ and $x$ or $w$ are not recolored by $c(y)$ then, its color remains unchanged. If $c(y)\leq p-2$ and $x$ is recolored by $c(y)$ then, we assign the color $p$ to $y$. If $c(y) \leq p-2$ and $w$ is recolored by $c(y)$ then, we recolor $y$ by $p-1$. If none of above situations happens then, assign the greatest color in $S$ to $y$ as its new color. Update $S$ and go to another vertex in $L_i$. Note that each recolored vertex $u$ of $L_{i-1}$ with $c(u)\geq p-1$ has all colors $1, \ldots, p$ in its closed neighborhood.

\noindent We continue the recoloring process until we confront with at least $p$ vertices whose primary colors (i.e. in $C$) is at least $p$. It will be shown, later in a claim, that if we continue the recoloring process until at most the level $\lfloor \log k \rfloor +3$, i.e. $L_{\lfloor \log k \rfloor +3}$, then $p$ vertices of color $\geq p$ are surely visited.

\noindent {\bf Pruning stage:}

\noindent In this stage we prune (remove) some vertices from the recolored subgraph.

\noindent Note that the desired induced subgraph we are looking for, is finally constructed on the remained vertices. Therefore the following removal of some vertices does not make any problem. Starting from $L_3$, we remove all vertices $w$ in $L_3$ satisfying the following properties. For some $i \leq p-2$, $w \in CH(v_i)$ and $w \notin CH(v_j)$, for any $j > p-2$. In general, we remove all vertices $w$ in $L_i$, $i\geq 4$, satisfying the following properties. For some $y$ with $c(y) \leq p-2$, $w \in CH(y)$ but $w \notin CH(y')$, for any $y'$ with $p-2 < c(y')$. The pruning is continued until the level $\lfloor \log k \rfloor +3$. Hence, for each $i\in \{3, 4, \ldots, \lfloor \log k \rfloor +3\}$ the level $L_i$ is pruned. Denote by $L'_i$ the remaining vertices in $L_i$.

\noindent \textbf{Properness of the recoloring:}

\noindent As we mentioned before some adjacent vertices may receive the same color in the recoloring process. We explain now some additional local recoloring to resolve this problem. Recall that the primary color of each vertex is denoted by $c(v)$ and in case that its color is changed then its new color is denoted by $c'(v)$. The recoloring of $L_2$ remains proper unless for some $i$ and $j$, $v_i$ is adjacent to $v_j$ and $c(v_j)\not=c'(v_j)=i=c(v_i)=c'(v_i)$. Note that this is the only case where properness of color $i$ in $L_2$ fails, because each color $c \in S$ is assigned to just one vertex of $v_i$ with primary color $i \geq p-1$, and hence no pair of vertices from $v_{p-1}, \ldots, v_{k-1}$ receive the same color from $S$. In this situation, it's enough to remove $v_i$ from $L'_2$.

\noindent Assume that the recolored $C'$ is proper on $L'_1\cup L'_2 \cup \ldots \cup L'_{i-1}$. Now we discuss the vertices in $L'_i$. The recoloring is proper in this level too unless some vertex $v\in L'_i$ is adjacent to a vertex $u\in L'_{i-1}\cup L'_i$ with the same color, i.e. $c'(u)=c'(v)$. We apply the following final trick in order to make the coloring proper.

\noindent \textbf{A trick to make the recoloring proper:}

\noindent Let $v \in CH(w)$ and $u \in L'_{j}$,  $j \in \{i,i-1\}$. We have the following possibilities.

\noindent Case 1. $u \notin CH(w)$.
\begin{enumerate}

\item[\textbf{i.}]
If $c(v) \geq {p-1}$. In this case we have either $c(v) \neq c'(v) = c'(u)=c(u) \leq p-2$ or $c(v) =c'(v) = c'(u) \neq c(u)$.
In the first case, $c(v) \neq c'(v) = c'(u)=c(u) \leq p-2$, let $u\in CH(x)$. Since we have the vertex $u$ in this case, we know that the vertices of $CH(x)$ are not removed in the pruning stage. Therefore $c(x) \geq p-1$ and consequently we have $|CH(x)| \geq p-2$. For any $k\geq 12$, we have $p-2 \geq 4$, which means there are at least $4$ vertices in $CH(x)$ with primary color less than $p-1$. On the other hand, based on the recoloring method, at most two vertices among $CH(x)$ with primary color less than $p-1$ are recolored. It follows that there exists a vertex $y$ such that  $y \in CH(x)$, $y\neq u$ and $c'(y) = c(y)\leq p-2$. Also note that by Lemma \ref{H-cactus}, $y$ is not adjacent to any vertex of $L'_{i-1}\cup L'_i$ (except $x$). In this case exchange the colors of $u$ and $y$. In case that $k\in \{10, 11\}$, we have already shown that the recoloring is finished in $L_2$. Hence such a possibility does not happen.

\noindent In the case $c(v) =c'(v) = c'(u) \neq c(u)$, take a $z \in CH(w)$ with $c(z)\leq p-2$. Such a vertex exists because $k>9$. Also by Lemma \ref{H-cactus}, $z$ is not adjacent to any vertex of $ L'_{i-1}\cup L'_i$ (except $w$). Exchange the colors of $v$ and $z$.
\item[\textbf{ii.}]
If $c(v) < p-1$.  Let $x$ be any vertex in $CH(w)$ with $c(x)\leq p-2$. Such a vertex
exists because $k>9$. Also by Lemma \ref{H-cactus}, $x$ is not adjacent to any vertex of $ L'_{i-1}\cup L'_i$ (except $w$). Exchange the colors of $v$ and $x$.
\end{enumerate}

\noindent Case 2. $u \in CH(w)$. In this case $v$ and $u$ are two adjacent vertices in $CH(w)$ and we may assume that $c'(u)=c(u)=c'(v)\not=c(v)$. Remove $u$ from $L'_i$.

\noindent Case 3. $u=w$. This case happens only when $c(v)=c'(v)$, $c(u) \not=c'(u)$. In this case remove $v$ from $L'_i$.

\noindent Denote by $D$ the set of all vertices $w$ with $c(w)\geq p-1$ and $c(w)\not= c'(w)$. Denote by $G'$ the subgraph of $G$ induced on the remaining vertices in $L'_1 \cup L'_2 \cup \ldots \cup L'_{\lfloor \log{k} \rfloor+3}$. Note that for each vertex $w\in D$ any color from $\{1, 2, \ldots , p\}$ appears in its closed neighborhood $N[w]$ in $G'$. Also note that by the method of recoloring, when a vertex receives a color of $S$, the color is removed from $S$ and  therefore no other vertex gets the same color in future. In other words, there are not two vertices in $D$ with a same color.

\noindent To complete the proof, it suffices to show $|D|\geq p$. By proving this fact we conclude that $G'$ admits a {\rm b}-coloring using $p$ colors, as desired.

\noindent \textbf{Claim:} $\left| D \right| \geq p$.

\noindent \textbf{Proof of the claim:}
We count the vertices in $L'_1 \cup L'_2 \cup \ldots \cup L'_{\lfloor \log{k} \rfloor+3}$ of primary color $\geq p-1$. To do so, we first show by the induction on $i\geq 2$ that for any primary color $j\geq p-1$, there are at least $k-j-i+1 \choose i-2$ many vertices in  $L'_{i} \cap D$. We first show validity of the induction hypothesis for $i=2$.

\noindent Regarding to the recoloring of vertices in $L'_2$, note that $v_{p-1} , \ldots , v_{k-1}$ satisfy membership in $D$. In other words, there are exactly one vertex of each color $p-1 \leq j \leq k-1$ in $L'_2 \cap D$. This proves the assertion for $i=2$.
	
\noindent Suppose now that the induction hypothesis holds in all levels $L'_1, \ldots, L'_i$. We prove it for the next level $L'_{i+1}$. We have $k-j-i+1 \choose i-2$ vertices of primary color $j$ in $L'_i \cap D$. We have to show that there are  $k-j-i \choose i-1$ vertices of primary color $j$ in $L'_{i+1} \cap D$.

\noindent First we are going to explain that there are at least  $k-p-(2i-5)$ distinct colors in $L'_i \cap D$ all greater than $p-1$ and the largest one is at least $k-2i+3$.
In each level, a vertex with a greater color (a greater number in fact), requires more number of neighbors. This forces the next level to contain more vertices with larger colors. And this situation will be repeated for the other levels. Therefore, since we desire to obtain a lower bound for $|D|$, we estimate the possible smallest size for $D$ among all possibilities for the edge connections in the subgraph of $G$ induced on $L'_1 \cup L'_2 \cup \ldots \cup L'_{\lfloor \log{k} \rfloor+3}$. Hence, in the following we look for a situation in which the maximum color in each level $L'_2, \ldots, L'_{\lfloor \log{k} \rfloor+3}$ is the possible smallest value, provided that the whole coloring satisfies the Grundy coloring properties.

\noindent Based on the method of recoloring, there are all colors $p-1, \ldots, k-1$ in $L'_2$. In  $L'_3$, $k-2$ is the greatest color which is possible to appear (as a child of $v_{k-1}$).
Because of the above-mentioned minimality requirement for the cardinality of $D$, we may assume that $v_{k-1}$ is adjacent to $v_{k-2}$ and therefore the color $k-2$ does not appear in $L'_3$. On the other hand, among the vertices of $L'_2$, no more than two vertices are adjacent (because the graph is a cactus). It means that the colors $p-1, \ldots, k-3$ appear in $L'_3$ (actually we omitted the case in which color $k-2$ appears in $L'_i$).
Similarly, in an arbitrary level $L'_i$, assume that $c$ is the greatest possible color which can be existed in $L'_i$. Consider a case in which the vertex with color $c$ has its neighbor with color $c-1$ in the same level, $L'_i$.
This way, we describe a coloring in which the greatest color of a level is the smallest possible one. In fact we describe the aforesaid situation in which $D$ has the smallest possible size. Doing this for $L'_2, \ldots, L'_i$, the greatest color in $L'_i$ is at least $k-2i+3$. Therefore there are at least $(k-2i+3)-(p-2)=k-p-2i+5$ distinct colors from color $p-1$ to the color $k-2i+3$ in $L'_i$.

\noindent In $L'_i$ the vertices are the children sets (i.e. all children of a vertex) of the vertices in $L'_{i-1}$. Among the vertices of a children set of color at least $j+1$, at most one vertex is adjacent to a color $j$ in $L'_1 \cup L'_2 \cup \ldots \cup L'_{i}$ (by Lemma \ref{H-cactus}). Therefore in a children set, at most one vertex with color at least $j+1$ does not introduce a color $j$ to $L'_{i+1}$. Assume that in every children set there is a vertex which does not introduce color $j$ to $L'_{i+1}$. Without loss of generality, assume that none of the vertices with color $j+1$ introduce a color $j$ to $L'_{i+1}$ (in fact this is a situation with the smallest number of colors $j$ in $L'_{i+1}$). Hence each vertex of $L'_i$ of color $c$, $c \in \{j+2, \ldots, k-2i+3\}$, introduces one vertex of color $j$ to be put in $L'_{i+1} \cap D$. It follows that the number of vertices of color $j$ in $L'_{i+1} \cap D$ is at least

\begin{center}
 $\underbrace{{k-j-i-1 \choose i-2}}_\text{ vertices of color $j+2$} + \underbrace{{k-j-i-2 \choose i-2}}_\text{ vertices of color $j+3$} + \cdots + \underbrace{{i-2 \choose i-2}}_\text{ vertices of color $k-2i+3$}= {k-j-i \choose i-1}$.
\end{center}

\noindent In the following we count the total number of vertices in a level $L'_i \cap D$.
\noindent For every $i$, there are at least $k-j-i+1 \choose i-2$ many vertices in $L'_{i} \cap D$ with primary color $j$, $p-1\leq j$. Note that there are at least  $k-p-(2i-5)$ distinct colors in $L'_{i} \cap D$ all more than $p-1$ and hence the greatest one is at least $k-2i+3$.
This argument shows that there are overall at least
\begin{center}
$\sum_{j=p-1}^{k-2i+3}{k-j-i+1 \choose i-2}={k-p-i+3 \choose i-1}$
\end{center}
vertices in $L'_i \cap D$. We end this process in $L'_{\lfloor \log k \rfloor +3}$. Note that none of the removed vertices in the pruning process, belongs to $D$, so removal of them does not decrease the cardinality of $D$.
	
\noindent We obtain the following lower bound for $D$
 \begin{center}
$|D|\geq \sum_{i=1}^{\lfloor \log k \rfloor +3}{k-p-i+3 \choose i-1}$.
\end{center}
\noindent By replacing $k=2^m +t, 0 \leq t \leq 2^m-1 $, and  $ k-p=2m$ we have
\begin{center}
$|D|\geq \sum_{i=1}^{m+3}{2m-i+3 \choose i-1}\geq \sum_{i=1}^{m+3}{m+1 \choose i-1}=2^{m+1}\geq 2^{m}+t-2m = k- 2m=p$.
\end{center}

\noindent It implies that $G'$ admits a {\rm b}-coloring using $p$ colors, therefore ${\rm b}(G') \geq k-2 \lfloor \log{k} \rfloor$ and so ${\rm{b}}(G) \geq k- 2\lfloor \log{k} \rfloor$.
\end{proof}

\noindent In the following we present an example of a $b$-monotone cactus graph having the Grundy  number $k$ and arbitrarily large ${\rm b}$-chromatic number.

\noindent For two positive integers $4\leq k$ and $t$, $k< t$, consider a set of $t$ vertices $v_1, \ldots, v_t$ such that every pair of  vertices $v_i$ and $v_{i+1}$, $1 \leq i \leq t-1$, are adjacent and there is no other edge among the vertices. Attach $t-1$ vertices $u_{i,j}, 1 \leq j \leq t \; j \neq i$, to a vertex $v_i$ for all $1 \leq i \leq t$. Also attach a $k$-atom $T_k$  by amalgamating $u_{2,1}$ with one vertex from $T_k$ which has degree one. Call this graph $G$.
There are not $t+1$ vertices with degree at least $t$ in $G$. Therefore ${\rm{b}}(G)\leq t$. Now we present a $\rm{b}$-coloring for $G$ with $t$ colors.  Assign colors $1, \ldots, t$, respectively to the vertices $v_1, \ldots, v_t$. Then for every vertex $u_{i,j}$ put a color $j$. Now with any proper $k$-coloring of the subgraph $T_k$, in which $u_{2,1}$ has color $1$, we have a $\rm{b}$-coloring of $G$ with $t$ colors.
On the other hand, for a Grundy coloring of $G$, the vertices $u_{i,j}$ (except $u_{2,1}$) can not have color more than 2
and therefore the vertices $v_i,  1\leq i \leq t$, can not have color more than 5. Also no vertex in the subgraph $T_k$ can be colored with a color greater than $k$. Therefore the Grundy number of $G$ is at most $k$.

\section{Some forbidden subgraphs}

\noindent Let $H_1, \ldots, H_k$ be any fixed set of graphs. A graph $G$ is $(H_1, \ldots, H_k)$-free if $G$ does not contain any induced subgraph isomorphic to $H_i$, for each $i$.
In the following by $K_4\setminus e$ we mean a graph obtained by removing any edge $e$ from $K_4$. Also $C_4$ stands for the cycle on four vertices. The next theorem proves Conjecture \ref{conj} for the family of $(K_4\setminus e, C_4)$-free graphs.

\begin{thm}
Let $G$ be a $(K_4\setminus e, C_4)$-free $b$-monotone graph. Then
$$\lfloor \Gamma(G)/2 \rfloor \leq {\rm b(G)}.$$\label{K4e}
\end{thm}

\noindent \begin{proof} Let $C$ be a Grundy coloring of $G$ using $k=\Gamma(G)$ colors and $H$ be the subgraph corresponding to $C$. Recall that $H$ consists of the levels $L_1, \ldots, L_t$ for some $t$ in which $L_1=\{v_k\}$ and $L_2=\{v_1, \ldots, v_{k-1}\}$. The color of $v_i$ is $i$ and $v_i$ is adjacent to $v_k$, for each $i\not=k$. Set $p= \lfloor k / 2 \rfloor$. Suppose first that $k$ is even.
In the following we recolor some vertices and at the same time remove some vertices in order to obtain a subgraph of $G$ with {\rm b}-chromatic number at least $p$.
In $L_1$, recolor $v_k$ by $p$. In $L_2$, assign colors $1, \ldots, p-1$ to $v_{p+1}, \ldots, v_{k-1}$, respectively and remove vertices $v_1, \ldots, v_p$. Note that for any $j$ with $1\leq j \leq k-p-1$, $v_{p+j}$ receives new color $j$. In other words, $c'(v_{p+j})=j$.

\noindent In general case, let $y \in CH(v_{p+l})$, $l \leq p-1$, be a vertex of $L_3$ with $c(y)=n$. If $l\not=n \leq p-1$ and the vertex $v_{p+l}$ is not adjacent to $v_{p+n}$ then keep $y$. Otherwise, remove $y$ from $L_3$. Note that $v_{p+n}$ has the new color $n$. Hence $v_{p+l}$ has all colors $1, \ldots, p$ in its closed neighborhood. Do this operation for all vertices of $CH(v_{p+l})$. We continue this process until $v_{k-1}$ and repeat the same technique for this vertex. We obtain that all the vertices $v_{p+1}, \ldots, v_{k-1}, v_k$ have colors $1, \ldots, p$ in their closed neighborhood. Denote by $L'_1, L'_2$ and $L'_3$, the remaining vertices from the levels $L_1, L_2$ and $L_3$, respectively. Define $G'=G[L'_1\cup L'_2 \cup L'_3]$.

\begin{figure}
	\begin{center}
		\begin{tikzpicture}
		\draw[black, thick] (0,0)-- (-1,-1);
		\draw[black, thick] (0,0)-- (1,-1);
		\draw[black, thick] (-1,-1)-- (0,-2);
		\draw[black, thick] (0,-2)-- (1,-1);
		
		\filldraw [black] (0,0) circle(2pt)
		node [anchor=south]{$ v_{k} $};
		\filldraw [black] (-1,-1) circle(2pt) node [anchor=east]{$v_{p+i}$};
		
		\filldraw [black] (1,-1) circle(2pt)
		node [anchor=west]{$w$};
		\filldraw [black] (0,-2) circle(2pt)
		node [anchor=north]{$z$};
		\end{tikzpicture}
		or
		\begin{tikzpicture}
		\draw[black, thick] (0,0)-- (-1,-1);
		\draw[black, thick] (0,0)-- (1,-1);
		\draw[black, thick] (-1,-1)-- (0,-2);
		\draw[black, thick] (0,-2)-- (1,-1);
		\draw[black, thick] (0,-2)-- (0,0);
		
		\filldraw [black] (0,0) circle(2pt)
		node [anchor=south]{$ v_{k} $};
		\filldraw [black] (-1,-1) circle(2pt) node [anchor=east]{$v_{p+i}$};
		
		\filldraw [black] (1,-1) circle(2pt)
		node [anchor=west]{$w$};
		\filldraw [black] (0,-2) circle(2pt)
		node [anchor=north]{$z$};
		\end{tikzpicture}
	\end{center}\label{picK4e}
	\caption{A situation in the proof of Theorem \ref{K4e}}
\end{figure}
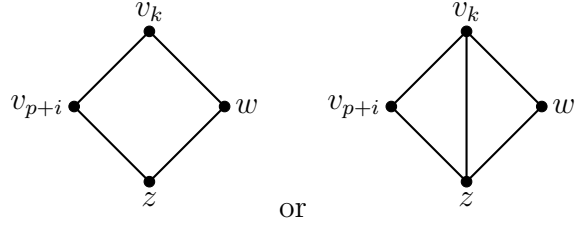

\noindent We claim that the resulting coloring for $G'$ is proper. The color of all vertices in $L'_1 \cup L'_2$ is changed and no color is assigned to more than one vertex in $L'_1 \cup L'_2$; therefore properness is met in these levels. Assume on the contrary that there is a vertex $z \in  L'_3$, with $c'(z)=i$, which is adjacent to a vertex $y$ with $c'(y)=i$. Since the color of vertices in $L'_3$ are not changed, then $c'(z)=c(z)$ and $y$ must be in $L'_2$. In fact $y=v_{p+i}$. According to the rules of the recoloring, there is a vertex $w \in L'_2$ which is not adjacent to $v_{p+i}$ and also $z \in CH(w)$. Here we obtain either an induced $C_4$ or an induced $K_4\setminus e$ on $\{v_k, v_{p+i}, z, w\}$, a contradiction. This situation is illustrated in Figure \ref{picK4e}. It follows that the recoloring is proper and $G'$ has a {\rm b}-coloring with $p$ colors and $\{v_k, v_{k-1}, \cdots, v_{p+1}\}$ as its color-dominating set. Since $G$ is $b$-monotone, we have $p\leq {\rm b}(G') \leq {\rm b}(G)$.

\noindent In the case that $k$ is an odd number, remove the class of vertices of color $k$ from $C$ and denote the resulting subgraph by $G_0$. We have $\Gamma(G_0)=k-1$ and $G_0$ satisfies the conditions of the previous case. We obtain $\lfloor \Gamma(G)/2 \rfloor = (k-1)/2 \leq {\rm b}(G_0) \leq {\rm b}(G)$.
\end{proof}

\noindent In the following we show that there is a $(K_4\setminus e, C_4)$-free graph $G$ such that ${\rm b}(G)=\lfloor \Gamma(G)/2 \rfloor +1$.

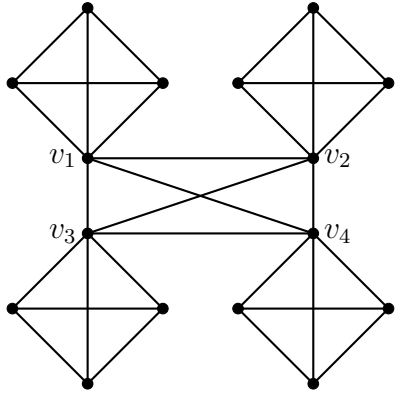
\begin{figure}
\begin{center}
\begin{tikzpicture}
\draw[black, thick] (0,2)-- (-1,1);
\draw[black, thick] (0,2)-- (1,1);
\draw[black, thick] (0,2)-- (0,0);
\draw[black, thick] (3,2)-- (2,1);
\draw[black, thick] (3,2)-- (4,1);
\draw[black, thick] (3,2)-- (3,0);
\draw[black, thick] (0,0)-- (0,-1);
\draw[black, thick] (0,0)-- (3,0);
\draw[black, thick] (0,0)-- (3,-1);
\draw[black, thick] (0,-1)-- (-1,-2);
\draw[black, thick] (0,-1)-- (1,-2);
\draw[black, thick] (0,-1)-- (0,-3);
\draw[black, thick] (3,-1)-- (2,-2);
\draw[black, thick] (3,-1)-- (3,-3);
\draw[black, thick] (3,-1)-- (4,-2);
\draw[black, thick] (3,-1)-- (3,0);
\draw[black, thick] (3,-1)-- (0,-1);
\draw[black, thick] (3,0)-- (0,-1);
\draw[black, thick] (4,-2)-- (2,-2);
\draw[black, thick] (3,-3)-- (2,-2);
\draw[black, thick] (3,-3)-- (4,-2);
\draw[black, thick] (4,1)-- (2,1);
\draw[black, thick] (1,-2)-- (-1,-2);
\draw[black, thick] (0,-3)-- (1,-2);
\draw[black, thick] (0,-3)-- (-1,-2);
\draw[black, thick] (-1,1)-- (1,1);
\draw[black, thick] (0,0)-- (-1,1);
\draw[black, thick] (0,0)-- (1,1);
\draw[black, thick] (3,0)-- (2,1);
\draw[black, thick] (3,0)-- (4,1);
\filldraw [black] (3,0) circle(2pt)
node [anchor=west]{$v_2$};

\filldraw [black] (0,0) circle(2pt)
node [anchor=east]{$v_1$};

\filldraw [black] (0,-1) circle(2pt)
node [anchor=east]{$v_3$};

\filldraw [black] (3,-1) circle(2pt)
node [anchor=west]{$v_4$};

\filldraw [black] (-1,1) circle(2pt)
node [anchor=west]{};
\filldraw [black] (2,-2) circle(2pt)
node [anchor=west]{};
\filldraw [black] (4,-2) circle(2pt)
node [anchor=west]{};
\filldraw [black] (3,-3) circle(2pt)
node [anchor=west]{};
\filldraw [black] (0,-3) circle(2pt)
node [anchor=west]{};
\filldraw [black] (1,-2) circle(2pt)
node [anchor=west]{};
\filldraw [black] (-1,-2) circle(2pt)
node [anchor=west]{};
\filldraw [black] (4,1) circle(2pt)
node [anchor=west]{};
\filldraw [black] (3,2) circle(2pt)
node [anchor=west]{};
\filldraw [black] (2,1) circle(2pt)
node [anchor=west]{};
\filldraw [black] (1,1) circle(2pt)
node [anchor=west]{};
\filldraw [black] (0,2) circle(2pt)
node [anchor=west]{};
\end{tikzpicture}
\end{center}
\caption{An (almost) sharpness example for Theorem \ref{K4e}}\label{sharpexam}
\end{figure}

\noindent For any positive integer $t$ we construct a graph $G_t$ as following. Consider a complete graph on $t$ vertices $v_1, \ldots, v_t$. Corresponding to each $v_i$ attach a complete graph $K(v_i)$ on $t$ vertices, by amalgamating $v_i$ with one vertex from $K(v_i)$. We have $V(K(v_i))\cap V(K(v_j))=\varnothing$, for every $i, j$ with $1 \leq  i< j \leq t$. Denote the resulting graph by $G_t$. The graph $G_4$ is depicted in Figure \ref{sharpexam}.
We show that $\Gamma(G_t)=2t-1$ and ${\rm b}(G_t)=t$.

\noindent We know that for any graph $G$, $\Gamma(G) \leq \Delta(G)+1$. $ \Delta(G_t)=2t-2$ hence $\Gamma(G_t)\leq 2t-1$. Now we present a Grundy coloring for $G_t$ with $2t-1$ colors. Assign colors $t, \ldots, 2t-1$ respectively to the vertices $v_1, \ldots, v_t$. Then for every vertex in $K(v_i), i\in \{1,\ldots t\}$, except $v_i$, assign a color $c$ from the set $\{1, \ldots, t-1\}$, such that no pair of vertices in $K(v_i)$ have the same color. This is a Grundy coloring of $G_t$ with $ 2t-1$ colors.

\noindent Since $m(G_t)=t$ and ${\rm{b}}(G_t)\leq m(G_t)$, we have ${\rm{b}}(G_t)\leq t$.  Now we present a $\rm{b}$-coloring for $G_t$ with $t$ colors.  For each $i \in \{1, \ldots, t \}$, assign color $ i $ to the vertex $v_i$. Then for every other vertex in $K(v_i), i\in \{1,\ldots t\}$, assign a color $c$ from the set $\{1, \ldots, t\} $, such that $c \neq i$. Also, it is easily seen that $G_t$ is $(K_4\setminus e, C_4)$-free.

\noindent The following corollary is a direct consequence of Theorem \ref{K4e} and the fact that any graph of girth at least $5$ is $b$-monotone.

\noindent \begin{cor}
Let $G$ be a graph with girth at least $5$. Then $\lfloor \Gamma(G)/2 \rfloor \leq {\rm b(G)}$.
\end{cor}

\noindent In the following theorem we obtain a result for graphs of girth at least $6$.

\begin{thm}\label{girth6}
Let $G$ be a graph of girth at least $6$. Then $\lfloor 2\Gamma(G)/3 \rfloor \leq {\rm b}(G)$.
\end{thm}

\noindent \begin{proof}
Set $\Gamma(G)=k$. Since every graph of girth at least $5$ is $b$-monotone then it suffices to find a subgraph $G'$ of $G$ with ${\rm b}(G')\geq \lfloor 2k/3 \rfloor$.
Let $C$ be a Grundy coloring of $G$ using $k$ colors and $H$ be the subgraph corresponding to $C$. Recall that $H$ is consisted of the levels $L_1, \ldots, L_t$ for some $t$ in which $L_1=\{v_k\}$ and $L_2=\{v_1, \ldots, v_{k-1}\}$. There are no edges between vertices of $L_2$, because of the girth. The color of $v_i$ is $i$ and $v_i$ is adjacent to $v_k$ for each $i\not=k$. Set $p= \lfloor 2k/3\rfloor$. Suppose first that $k=3t$ or $k=3t+2$, for some integer $t$.

\noindent Recolor $v_k$ by $p$ and $v_p, \ldots, v_{k-1}$ by $\lfloor k/3\rfloor, \ldots, p-1$, respectively. The color of other vertices in $L_2$ remain unchanged. To obtain $G'$ we put some vertices in $G'$ and ignore the rest of vertices. Let $w_{i,j}$ be an arbitrary vertex of color $j$ in $L_3$ and  $w_{i,j}\in CH(v_i)$. If $i<p$ then ignore $w_{i,j}$. Also for $p \leq i \leq k-2$ if $p \leq j$ or $j=c'(v_i)$, then ignore $w_{i,j}$. Otherwise, put $w_{i,j}$ in $G'$.
In the case that $i=k-1$, recolor all vertices $w_{k-1, \lceil 2k/3\rceil}, \ldots, w_{k-1, k-2}$ by $1, \ldots, \lfloor k/3\rfloor-1$, respectively and put them in $G'$. We do not change the colors of $w_{k-1,\lfloor \frac{k}{3}\rfloor}, \ldots, w_{k-1,p-2}$ and put these vertices too in $G'$. Then ignore all other vertices of $L_3$. Note that at this step, the resulting recoloring has the following property. For each $i$ with $p\leq i \leq k-1$, the vertex $v_i$  has all colors $1, \ldots, p$ in its closed neighborhood $N[v_i]$.

\noindent In $L_4$ for every $i \not= k-1$ remove $CH(w_{i,j})$. For $i=k-1$, let $x \in CH(w_{{k-1},j})$. For $j\in \{\lceil 2k/3\rceil, \ldots, k-2 \}$ keep $x$ unchanged when $c(x) \in \{1, \ldots, p\}$ and $c(x) \not= c'(w_{{k-1},j})$. Otherwise, remove $x$. Notice that if $k$ is a multiple of 3, we keep the child of $w_{k-1, \lceil 2k/3\rceil}$ of color $1$ and assign the new color $p$ to it.

\noindent Because the recoloring is occurred only in the first three levels, and since the girth is at least $6$, the resulting coloring is proper. Define  $D=\{ v_p, \ldots, v_{k-1}, v_{k} \} \cup \{w_{k-1,\lceil \frac{2k}{3}\rceil}, \ldots, w_{k-1,k-2}\}$. Note that the vertices of $D$ in subgraph $G'$ have all colors $\lbrace 1, 2, \ldots , p \rbrace$ on their closed neighborhood. It follows that $D$ is a color-dominating set for $G'$. Hence $G'$ is an induced subgraph in $G$ with ${\rm b}$-chromatic number at least $p$.

\noindent In the case that $k=3t+1$, for some integer $t$, remove the class of vertices of color $k$ in $C$ and denote by $G_0$ the resulting subgraph. We have $\Gamma(G_0)=k-1=3t$ and $G_0$ satisfies the conditions of the previous case. We obtain $\lfloor 2(3t)/3 \rfloor \leq {\rm b}(G_0)\leq {\rm b}(G)$. Then $\lfloor 2k/3 \rfloor=2t \leq {\rm b}(G)$.
\end{proof}

\noindent The result of Theorem \ref{girth6} can be generalized for higher girths. It can be proved using the same recoloring process that if $G$ has girth at least $8$ then $\lfloor 3\Gamma(G)/4 \rfloor \leq {\rm b}(G)$. We omit the details of the proof.

\section{Concluding Remarks}

\noindent One of the reviewers has mentioned that the inequality $\Gamma(G)\leq 2m(G)$ holds for all graphs $G$. Hence for any family of graphs in which $m(G)\leq {\rm b}(G)+1$, for any $G$ from the family, we immediately obtain $\Gamma(G)\leq 2{\rm b}(G)+2$. As confirmed by the reviewer, few families of graphs are known satisfying this property. Trees, cacti with $m\geq 7$ and graphs with girth at least 7, are some of such families, but we have very better comparative inequalities concerning the Grundy and the {\rm b}-chromatic number of graphs in these families. We present a proof of $\Gamma(G)\leq 2m(G)$ for any $G$. Let $C$ be a Grundy coloring of $G$ having $k=\Gamma(G)$ colors. Then there are vertices say $v_1, v_2, \ldots, v_k$ such that $d(v_1)\geq k-1, d(v_2)\geq k-2, \ldots, d(v_{k/2})\geq k/2$. Let $d_1\geq d_2 \geq \ldots \geq d_n$ be a vertex degree of $G$ in non-increasing form. It follows that $d_{k/2}\geq (k/2)-1$ and hence $m(G)\geq k/2$. Let ${\mathcal{F}}$ be any family such that for some function $f(.)$ and for any $G$ from the family, one has ${\rm b}(G)\geq f(m(G))$. It implies that Conjecture \ref{conj} is valid for ${\mathcal{F}}$.

\end{document}